\documentclass[letter, 10 pt, conference]{IEEEconf}

\usepackage{chaimath}
\usepackage{mathptmx} 
\usepackage{amsfonts}
\usepackage{times} 
\usepackage{amsmath} 
\usepackage{amssymb}  
\usepackage{amsthm}
\usepackage{amssymb, amsfonts, amsmath}
\usepackage{verbatim}
\usepackage{subfigure}
\usepackage{slashbox}
\usepackage{graphicx}

\title{{\bf A Sum-of-Squares Approach to the Analysis of Zeno Stability in Polynomial Hybrid Systems}}
\author{Chaitanya Murti and Matthew Peet%
\thanks{Chaitanya Murti is an M.S. student with the Cybernetic Systems and Control Lab (CSCL) and the Department of Electrical and Computer Engineering, Illinois Institute of Technology, 60616, USA, {\tt\footnotesize cmurti@hawk.iit.edu  }}%
\thanks{Matthew M. Peet is an Assistant Professor with the School of Engineering of Matter, Transport, and Energy, Arizona State University, Tempe, AZ, 85821 USA, {\tt\footnotesize mpeet@asu.edu  }}%
}
\date{\today}

\begin{document}
\maketitle

\begin{abstract}
Hybrid dynamical systems can exhibit many unique phenomena, such as Zeno behavior. Zeno behavior is the occurrence of infinite discrete transitions in finite time. Zeno behavior has been likened to a form of finite-time asymptotic stability, and corresponding Lyapunov theorems have been developed. In this paper, we propose a method to construct Lyapunov functions to prove Zeno stability of compact sets in cyclic hybrid systems with parametric uncertainties in the vector fields, domains and guard sets, and reset maps utilizing sum-of-squares programming. This technique can easily be applied to cyclic hybrid systems without parametric uncertainties as well. Examples illustrating the use of the proposed technique are also provided. 
\end{abstract}

\section{Introduction}
\noindent  Hybrid systems are dynamical systems with trajectories that exhibit both continuous flows and discrete transitions. As such, a variety of man-made systems can be modeled using the hybrid systems framework. Some exampleks are electrical systems with switching \cite{CTA521}, communication networks \cite{HespanhaHybCom}, embedded systems \cite{HybEmbed}, and air traffic control \cite{TomlinAir}. \\
\indent Recent research into hybrid systems has yielded results on stability of equilibria \cite{BranStability} and observability and controllability \cite{BemporadObsCont}. Several Lyapunov-based techniques for the analysis of hybrid systems, including the use of multiple Lyapunov functions \cite{BranMultLyap}, the construction of piecewise-quadratic Lyapunov functions \cite{JhnsQuadLyap}, and the utilization of Lyapunov techniques for robust stability analysis \cite{PettLennStab} have also been presented. More recently, a means to assess stability of hybrid systems by constructing higher-order polynomial Lyapunov functions using sum-of-squares techniques was presented in \cite{PrajPapach}, and a method to perform robust stability analysis using sum of squares techniques was provided in \cite{PapachPrajRobust}. However, there are still behaviors of hybrid systems that require further study. Among these phenomena are chattering and zeno behavior.  \\
\indent Zeno behavior is the occurrence of infinite transitions between discrete states in a finite period of time. Trajectories exhibiting this behavior are called Zeno executions, and converge to a set of points known as a Zeno equilibrium. Hybrid systems exhibiting Zeno behavior are described in detail in, for example, \cite{ZhangZunJohKarLygSasZeno}. Zeno behavior can cause simulations to halt or fail, since infinitely many transitions would need to be simulated, as noted in, e.g., \cite{ZhangZunJohKarLygSasZeno}. This problem was addressed in \cite{JohanssonReg} and \cite{LifeAftZenoAmes}, which describe methods to regularize hybrid systems to ensure that trajectories continue after the Zeno equilibrium. Sufficient conditions for Zeno behavior in first quadrant hybrid systems were given in \cite{AmesDFQ}, and further sufficient conditions for systems with nonlinear vector fields based on constant approximations were given in \cite{AmesNonlinconst}. More recently, necessary and sufficient Lyapunov conditions for the existence of isolated Zeno equilibria were first given in \cite{AmesLyap}. These results were extended in \cite{TeelGoebLyap}, where the concept of Zeno stability was described as an extension of finite-time asymptotic stability. Moreover, \cite{TeelGoebLyap} provided Lyapunov conditions for Zeno stability of compact sets. The results in \cite{TeelGoebLyap} were also shown to be equivalent to the theorem presented in \cite{AmesLyap}. The results of \cite{AmesLyap} were also exended by Ames and Lamperski to non-isolated Zeno equilibria in \cite{AmesLampTAC}. In a similar vein, a Lyapunov characterization of Filippov solutions was provided in \cite{AhmadiFilippov}. \\
\indent  In this paper, we use sum-of-squares programming to construct  Lyapunov functions which prove Zeno stability of compact sets, based on the results of \cite{AmesLyap} and \cite{TeelGoebLyap}. Moreover, the method presented in this paper allows for the verification of Zeno stability for polynomial hybrid systems with nonlinear polynomial vector fields and transitions. We also present a method to verify Zeno stability for systems with parametric uncertainties.\\
\indent The outline of the paper is as follows: in Section II, definitions of sum-of-squares polynomials, hybrid systems and their executions, and Zeno executions and equilibria are presented. Section II also details Lyapunov conditions for Zeno stability as described in \cite{AmesLyap}. In Section III, we present a method to construct Lyapunov functions to prove Zeno stability of hybrid systems with parametric uncertainties using Sum-of-Squares optimization, and in Section IV, illustrative examples are provided.

\section{Preliminaries}
In this section, we provide a brief introduction to Sum-of-Squares polynomials and definitions for hybrid systems, their executions, and Zeno behavior.

\subsection{Sum of Squares Polynomials}
We use $\mathbf{R}[x]$ to denote the ring of polynomials generated by variables $x = (x_1,..,x_n)$.
\defnn{Sum of Squares Polynomial}{A polynomial $p(x):\Realsn{n}\rightarrow \Reals$ is said to be Sum of Squares (SOS) if there exist polynomials $f_i(x):\Realsn{n}\rightarrow \Reals $ such that
\[p(x) = \sum_i (f_i(x))^2\]
We use $p\in \Sigma_x\subset \mathbf{R}[x]$ to denote that $p$ is SOS.
}
The following result gives a polynomial-time complexity test to determine whether a polynomial is SOS.\par\vspace{-0.5cm}
\noindent\theoremC{\label{thmsos}
For a polynomial, $p$ of degree $2d$, $p\in \Sigma_x$ if and only if there exists a positive semidefinite matrix $Q$, such that
\[p(x) = Z(x)^TQZ(x)\]
where $Z(x)$ is the vector of monomials of degree $d$ or less
}
Therefore, checking whether a polynomial is SOS is equivalent to checking the existence of a positive-semidefinite matrix $Q$ under some affine constraints, which can be solved with semidefinite programming. Thus, while checking polynomial positivity is NP-hard, checking whether a polynomial is SOS is decidable in polynomial time.
\\
\indent In this paper, Positivstellensatz results from algebraic geometry are used extensively to create constraints that can be implemented using sum of squares programming. We use the Positivstellensatz to construct Lyapunov functions which are positive on bounded sets (see section IV).\par\noindent
For further details and proofs, we refer to \cite{Stengle} and \cite{ParilloThs}.
\subsection{Hybrid Systems}
In this section, we define hybrid systems and their executions. We use similar notation to that given in \cite{HybVDSSchu} and, more recently, \cite{AmesLyap}.
\defnn{Hybrid System}{
A hybrid system $H$ is a tuple $H = (Q,E,D,F,G,R)$ where
\begin{itemize}
\item $Q$ is a finite collection of discrete states or indices 
\item $E\subset Q\times Q$ is a collection of edges, where for any edge $e = (q,q')$ we use the functions $s$ and $t$ to denote the start and end, so that for $e = (q,q')$, $s(e) = q$ and $t(e) = q'$
\item  $D = \{D_q\}_{q\in Q}$ is a collection of Domains, where for each $q\in Q$, $D_q\subseteq \Realsn{n}$
\item $F = \{f_q\}_{q\in Q}$ is a collection of vector fields, where for each $q\in Q$, $f_q:D_q \rightarrow \Realsn{n}$
\item $G = \{G_{e}\}_{e\in E}$ is a collection of guard sets, where for each $e = (q,q') \in E$, $G_{e} \subset D_{q}$
\item $R = \{\phi_{e}\}_{e\in E}$ is a collection of Reset Maps,  where for each $e=(q,q')\in E$, $\phi_{e}:G_{e}\rightarrow D_{q'}$.
\end{itemize}
}

\defnn{Cyclic hybrid system}{A cyclic hybrid system $H_c$ is a hybrid system where for each domain $q \in Q$, we can associate a unique edge $e(q) = (q,q_i) \in E$ such that $s(e(q))=q$ and such that for any $q \in Q$, $q=t(e(t(e(\cdots t(e(t(e(q))))))))$. That is, the set of edges forms a directed graph.
}
 
\defnn{Hybrid System Execution}{Consider the tuple
$\chi = (I,T,p,C)$ where 
\begin{itemize}
\item $I\subseteq \mathbb{N}$ is index of intervals
\item $T = \{T_i\}_{i\in I}$ are a set of open time intervals associated with points in time $\tau_i$ as $T_i = (\tau_i,\tau_{i+1}) \subset \Realsn{n}_+$ where $T_{i+1} = (\tau_{i+1},\tau_{i+2})$
\item $p:I\rightarrow Q$ maps each interval to a domain,
\item $C = \{c_{i}(t)\}_{i\in I}$ is a set of continuously differentiable functions. We say $\chi$ is an \emph{execution} of the hybrid system $H=(Q,E,D,F,G,R)$ with initial condition $(q_0,x_0)$ if
 $c_{1}(0) = x_0$ and $p(1) = q_0$.
\item $\dot{c}_{i}(t)= f_{p(i)}(c_{i}(t))$ for $t\in T_i$ and for all $i \in I$; $c_i(t)\in D_{p(i)}$ for $t\in T_i$ and for all $i \in I$; $c_i(\tau_{i+1})\in G_{(p(i),p(i+1))}$ for all $i\in I$; $c_{i+1}(\tau_{t+1}) = \phi_{(p(i),p(i+1))}(c_i(\tau_i))$ for all $i\in I$.
\end{itemize}

}
\subsection{Zeno Stability in Hybrid Dynamical Systems}
We now present definitions of Zeno executions, equilibria, and stability, along with necessary and sufficient conditions for Zeno stability as presented in \cite{AmesLyap} and \cite{TeelGoebLyap}.
\defnn{Zeno Execution}{
We say an execution $\chi = (I,T,p,C)$ starting from $(q_0,x_0)$ of a hybrid System $=(Q,E,D,F,G,R)$ is Zeno if
\begin{enumerate}
\item $I = \mathbb{N}$
\item $\lim_{i \rightarrow \infty} \tau_i < \infty$
\end{enumerate}
}
\defnn{Zeno Equilibrium}{A set $z = \{z_q\}_{q\in Q}$ is a Zeno equilibrium of a Hybrid System $H=(Q,E,D,F,G,R)$ if it satisfies
\begin{enumerate}
\item For each edge $e = (q,q') \in E$, $z_q\in G_e$ and $\phi_e(z_q) = z_{q'}$.
\item $f_q(z_q)\neq 0$ for all $q \in Q$.
\end{enumerate}
 Note that for any $z\in \{z_q\}_{q\in Q}$, where $\{z_q\}_{q\in Q}$ is a Zeno equilibrium of a cyclic hybrid system $H_c$,
\[\left (\phi_{i-1}\circ \cdots \circ \phi_0 \cdots \phi_i\right ) (z) = z\]
}
Next, we define Zeno stability:
\defnn{Zeno Stability}{
Let $H = (Q,E,D,F,G,R)$ be a hybrid system, and let $z = \{z_q\}_{q\in Q}$ be a compact set.  The set $z$ is Zeno stable if, for each $q\in Q$, there exist neighborhoods $Z_q$, where $z_q\in Z_q$, such that for any initial condition $x_0\in \bigcup_{q\in Q}Z_q$, the execution $\chi = (I,T,p,C)$, with $c_o(t_0) = x_0$ is Zeno, and converges to $z$.
}\par
\noindent Note that this definition of Zeno stability is consistent with the stability definitions provided in \cite{TeelGoebLyap}.
We now reiterate the Lyapunov conditions for the stability of Zeno equilibria in cyclic hybrid systems presented in \cite{AmesLyap}, which are as follows:\par

\theoremnC{Lamperski and Ames}{\label{thm2}
\noindent Consider a cyclic hybrid system $H = (Q,E,D,F,G,R)$, with an isolated Zeno equilibrium $\{z_q\}_{q\in Q}$. Let $\{W_q\}_{q\in Q}$ be a collection of open neighborhoods of $\{z_q\}_{q\in Q}$. Suppose there exist continuously differentiable  functions $V_q:\Realsn{n}\rightarrow \Reals$ and $B_q:\Realsn{n}\rightarrow \Reals$, and non-negative constants $\{r_q\}_{q\in Q}$, $\gamma_a$, and $\gamma_b$, where $r_q \in [0,1]$, and $r_q<1$ for some $q$ and such that
\begin{align}
&V_q(x) &&> 0 \quad \text{for all } x\in W_q\backslash z_q, q \in Q\label{eq:EC1}\\
&V_q(z_q) &&= 0, \quad \text{for all } q \in Q\;\label{eq:EC1b}\\
&\grad {V}_q^T(x)f_q(x) &&\leq 0\quad \text{for all } x\in W_q, \, q \in Q\label{eq:EC2}\\
&B_q(x) &&\geq 0\quad \text{for all } x\in W_q, \, q \in Q \label{eq:EC3}\\
&\grad {B}_q^T(x)f_q(x) &&< 0\quad \text{for all } x\in W_q, \, q \in Q  \label{eq:EC4}\\
&V_{q'}(R_{(q,q')}(x)) &&\leq r_q V_q(x),\label{eq:C1}\\
&&&\quad \text{for all } e = (q,q') \in E \text{ and } x\in G_{e}\cap W_q \notag\\
&B_q(R_{(q',q)}(x)) &&\leq \gamma_b \left(V_q(R_{(q,q')}(x)) \right)^{\gamma_a}\label{eq:C2} \\
&&& \text{for all } e =(q,q') \in E \text{ and }x \in G_{e} \cap W_q.\notag \vspace{-8mm}
\end{align}
Then $\{z_q\}_{q\in Q}$ is Zeno stable.
}
\par
\noindent As noted in \cite{TeelGoebLyap}, the conditions above are equivalent to those given in \cite[Proposition 5.2]{TeelGoebLyap}. Thus, satisfying EC1-C2 is also sufficient to prove asymptotic Zeno stability of a compact set. This in turn allows us to relax the restriction $f_q(z_q)\neq 0$.\par 
 To simplify notation, we will use the sufficient conditions of Theorem 3 as follows. Note that our subsequent analysis can be easily applied directly to the conditions of Theorem 2 and in our numerical examples we have tested both sets of conditions and they yield similar results. \par
\theoremC{\label{thm2}
Let $H = (Q,E,D,F,G,R)$ be a cyclic hybrid system, and let $z = \{z_q\}_{q\in Q}$ be a compact set. Let $\{W_q\subset D_q\}_{q\in Q}$, be a collection of neighborhoods of the $\{z_q\}_{q\in Q}$. Suppose that there exist continuously differentiable functions $V_q:W_q\rightarrow \Reals$, and positive constants $\{r_q\}_{q\in Q}$ and $\gamma$, where $r_q \in (0,1]$, and $r_q<1$ for some $q$ and such that
\begin{align}
&V_q(x) &&> 0 \quad \text{for all } x\in W_q\backslash z_q, q \in Q \label{eq:NC1}\\
&V_q(z_q) &&= 0, \quad \text{for all } q \in Q \label{eq:NC2}\\
&\grad {V}_q^T(x)f_q(x) &&\leq-\gamma\;\quad \text{for all } x\in W_q, \, q \in Q \label{eq:NC3}\\
&r_q V_q(x) &&\geq V_{q'}(\phi_{e}(x)) \label{eq:NC4}\\
&&& \text{for all } e =(q,q') \in E \text{ and }x \in G_{e} \cap W_q.\notag
\end{align}
then $z$ is Zeno stable.
}
\ProofC{
We show that if for each $q\in Q$, we can find a $V_q$ such that (\ref{eq:NC1})-(\ref{eq:NC4}) are satisfied, then the same $V_q$ also satisfies (\ref{eq:EC1})-(\ref{eq:C2}). From inspection, it is clear that if $V_q$ satisfies (\ref{eq:NC1})-(\ref{eq:NC4}), then (\ref{eq:EC1})-(\ref{eq:EC2}) and (\ref{eq:C1}) are satisfied. Second, choose $B_q=V_q$ for each $q\in Q$. From inspection, it is clear that $V_q$ also satisfies (\ref{eq:EC3}) and (\ref{eq:EC4}). Last, if $\gamma_a = \gamma_b=1$, we get $V_q\leq V_q$, where the equality holds. From this, we see that for each $q\in Q$, $V_q$ also satisfies (\ref{eq:C2}). Thus, the theorem is proved.
}
\section{Using Sum-of-Squares Programming to prove Zeno Stability}
Theorem 3 provides sufficient conditions for Zeno stability in cyclic hybrid systems.  We now show that these conditions can be enforced using SOS, even for systems with parametric uncertainties. First, define the vector of parametric uncertainties $P$ to lie within a semialgebraic set 
\begin{equation}\label{eq:params}
P:=\{p\in\Reals: \tilde{p}_{k}(p)\geq 0, k=1,2,...,K_1 \}.
\end{equation}
We then present the following assumption:
\assn{
\noindent For the purposes of this paper, we consider hybrid systems with polynomial vector fields and resets, and semialgebraic domains and guard sets, with parametric uncertainties in the each of the above. Let $P$ be defined as in (\ref{eq:params}). We implicitly assume that associated with every hybrid system is a set of polynomials $g_{qi}(x,p)$, $h_{e,k}(x,p)$ for $q \in Q$, $e \in E$, $i=k=1,\cdots,K_{q}$ and $k = 1,\cdots,N_{q}$ for some $K_q,N_q>0$, and $p\in P$. \\ In this framework, the domains of the hybrid system $H$ are defined as
\begin{equation}
D_q = \{x\in \Realsn{n}\;: \;g_{qk}(x,p)\geq 0,\; k=1,2,\cdots,K_{q}\}
\end{equation}
where $g_{qk}\in \mathbf{R}[x,p]$, $K_{q}\in \mathbb{N}$, and $p\in P$. The guard sets are defined as
\begin{equation}
G_e = \{x\in \Realsn{n}\;:\; h_{e,0}(x,p) = 0,\; h_{e,k}(x,p)\geq 0,\;k = 1,2,\cdots,N_{q} \}
\end{equation}
where each $h_{ek}\in\mathbf{R}[x,p]$, $N_q\in \mathbb{N}$, and $p\in P$. Lastly, for each $e=(q,q')\in E$, the reset map $\phi_{e}$ is given by the vector-valued polynomial function
\begin{equation}
\phi_{e} = [\phi_{e,1}(x,p),\cdots,\phi_{e,n}(x,p)]^T
\end{equation}
where $\phi_{e,j} \in \mathbf{R}[x,p]$ for $j =1 ,\cdots,n$, and $p\in P$.}

Let $H = (Q,E,D,F,G,R)$ be a cyclic hybrid system, and let $z=\{z_q\}_{q\in Q}$ be a compact set. Let $\{W_q\}_{q\in Q}$ be a collection of neighborhoods of $\{z_q\}_{q\in Q}$. We consider $W_q$ of the form
\[W_q:=\{x\in\Realsn{n}: w_{qk}(x)>0, k=1,2,...,K_q\}\]
where each $w_{qk}(x) \in \mathbf{R}[x]$.\par
\noindent Consider Feasibility Problem 1:\par \noindent
\noindent{\bf Feasibility Problem 1:}\par
For hybrid system $H = (Q,E,D,F,G,R)$, find
\begin{itemize}
\item $a_{qk}$, $c_{qk}$, $i_{qk}$, $\in \Sigma_{x,p}$, for $k=1,2,...,K_{qw}$, $p\in P$ and $q \in Q$;
\item $b_{qk}$, $d_{qk}$, $j_{qk} \in \Sigma_{x,p}$, for $k=1,2,...,K_q$, $p\in P$, and $q \in Q$.
\item $ \eta_{qk}$, $ \beta_{qk}$, $ \zeta_{qk} \in \Sigma_{x,p}$, for $k=1,2,...,K_1$, $p\in P$, and $q \in Q$.
\item $m_{e,l} \in \Sigma_{x,p}$ for $e \in E$, $p\in P$, and $l=1,2,...,N_q$
\item $V_q$, $m_{e,0}\in \mathbf{R}[x,p]$ for $e \in E$, $p\in P$, and $q \in Q$.
\item Constants $\alpha, \gamma > 0$, $\{r_q\}_{q\in Q}\in (0,1]$ such that $r_q<1$ for some $q \in Q$.
\end{itemize}
such that
\begin{flalign}
\notag &V_q - \alpha x^T x-\sum_{k=1}^{K_{qw}} a_{qk}w_{qk} -\sum_{k=1}^{K_q} b_{qk}g_{qk}\\
 & \qquad \qquad \qquad - \sum_{k_1=1}^{K_1} \eta_{qk_1}\tilde{p}_{qk} \in \Sigma_{x,p} \quad \text{for all }q\in Q\label{eq:R1}\\
 &V_q(z_q,p) = 0 \quad \text{ for all }q \in Q \label{eq:R2}\\
\notag &-\grad V_q^Tf_q - \gamma -\sum_{k=1}^{K_{qw}} c_{qk}w_{qk} -\sum_{k=1}^{K_q} d_{qk}g_{qk}\\
 & \qquad \qquad \qquad - \sum_{k_1=1}^{K_1} \beta_{qk_1}\tilde{p}_{qk} \in \Sigma_{x,p} \quad \text{ for all } q \in Q\label{eq:R3}\\
\nonumber&r_q V_q - V_{q'}(\phi_{e})-  m_{e,0}h_{e,0} - \sum_{l=1}^{N_q} m_{e,l}h_{e,l}- \sum_{k=1}^{K_{qw}} i_{qk}w_{qk}\\
 & \qquad  - \sum_{k=1}^{K_{q}} j_{qk}g_{qk} - \sum_{k=1}^{K_1} \zeta_{qk}\tilde{p}_{qk} \in \Sigma_{x,p} \;\;  \text{ for all } e = (q,q') \in E. \label{eq:R4}
\end{flalign}

\par\noindent
\theoremC{\label{robust}
Consider a cyclic hybrid system $H = (Q,E,D,F,G,R)$, and let $z=\{z_q\}_{q\in Q}$ be a compact set. If Feasibility Problem 2 has a solution, then $z$ is Zeno stable for all $p\in P$.
}
\ProofC{

To prove the theorem we show that if $V_q$, $q \in Q$ are elements of a solution of Feasibility Problem 1, then for each $q\in Q$, the same $V_q$ also satisfy (\ref{eq:NC1})-(\ref{eq:NC4}) of Theorem 3 for all $p\in P$.
That is, we show that if the $V_q$ satisfy (\ref{eq:R1})-(\ref{eq:R4}), then the same $V_q$ also satisfies (\ref{eq:NC1})-(\ref{eq:NC4}) for all $p\in P$.  \\
First, we observe that~\eqref{eq:R2} directly implies~\eqref{eq:NC2} for $p\in P$. Next, from~\eqref{eq:R1}, we know that
\begin{align*}
& V_q(x,p) \geq \sum_{k=1}^{K_{qw}} a_{qk}(x,p)w_{qk}(x) 
+\sum_{k=1}^{K_q}b_{qk}(x,p)g_{qk}(x,p) + \alpha x^T x +\\& \qquad\quad+ \sum_{k_1=1}^{K_1} \eta_{qk_1}(x,p)\tilde{p}_{qk}(p)
\end{align*}
Since $a_{qk}(x,p)$, $b_{qk}(x,p)$, and $\eta_{qk}(x,p)$ are SOS, and thus, always nonnegative, by the Positivstellensatz and the definitions of $W_q$, $P$, and $D_q$, we have that $V_q(x) \geq \alpha x^T x$ for all $x\in W_q\subset D_q$ and all $p\in P$. Thus, (\ref{eq:R1}) implies (\ref{eq:NC1}) is satisfied. Similarly, from (\ref{eq:R4}),
\begin{align*}
 & -\grad V_q^T(x,p)f_q(x,p) -\gamma\geq \sum_{k=1}^{K_{qw}} c_{qk}(x,p)w_{qk}(x)+\\&\sum_{k=1}^{K_q} d_{qk}(x,p)g_{qk}(x,p) - \sum_{k_1=1}^{K_1} \beta_{qk_1}\tilde{p}_{qk}.
\end{align*}
Since $c_{qk}(x,p)$ and $d_{qk}(x,p)$ are always nonnegative, by the definition of $P$, $D_q$ and $W_q$, $\grad V_q(x) ^T f_q(x,p) \le -\gamma$ for $x\in\{x\in\Realsn{n}: g_{qk}(x,p)\geq 0,\;w_{qk}(x,p)\geq 0\} = D_q \cap W_q$ and $p\in P$, which implies (\ref{eq:NC3}) is satisfied. Next, from~\eqref{eq:R4} we have that for all $e = (q,q') \in Q$,
\begin{align*}
& r_qV_q(x,p) - V_{q'}(\phi_{e}(x,p),p) \geq m_{e,0}(x,p)h_{e,0}(x,p) \\&\quad + \sum_{l=1}^{N_q} m_{e,l}(x,p)h_{e,l}(x,p)  + \sum_{k=1}^{K_q} i_{qk}(x,p)w_{qk}(x) \\&\qquad+ \sum_{k=1}^{K_q} j_{qk}(x,p)g_{qk}(x,p) +  \sum_{k=1}^{K_1} \zeta_{qk}(x,p)\tilde{p}_{qk}(p).
\end{align*}
First note that $h_{e,0}(x) = 0$ and hence $m_{e,0}(x,p)h_{e,0}(x,p) = 0$ on $G_e$. Since $m_{e,l} \in \Sigma_x$, we have  $m_{e,l}(x)h_{e,l}(x)\ge 0$ on $G_e$. Similarly  $j_{qk}(x,p)g_{qk}(x,p) \ge 0$ on $D_q$ and $i_{qk}(x,p)w_{qk}(x) \ge 0$ on $W_q$. It follows that $r_q V_q(x) - V_{q'}(\phi_{e}(x)) \geq 0$ when $x\in G_e \cap W_q \cap D_q$ for all $p\in P$ $e = (q,q') \in E$. Thus, we have shown that (\ref{eq:R4}) implies (\ref{eq:NC4}). \par
\noindent  Thus we conclude that the solution elements $V_q$ of Feasibility Problem 1 satisfy the conditions~(\ref{eq:NC1})-(\ref{eq:NC4}) of Theorem 4. Thus by Theorem 3 we conclude Zeno stability of $z$ for all $p\in P$.
}
\noindent{\bf Remark:\;\;} For systems without parametric uncertainty, we simply take the set of uncertain parameters to be empty. Thus,  all elements of Feasibility Problem 1 become dependent only on $x$, and all $\tilde{p}_{qk},\;\eta_{qk},\;\beta_{qk}$, and $\zeta_{qk}$ are $0$. A similar theorem for Zeno stability of hybrid systems without parametric uncertainty is stated explicitly in \cite{my_thesis}.

\section{Examples}
In this section, we provide some examples that illustrate the application of the given technique. We demonstrate Zeno stability in hybrid systems with polynomial vector fields and semialgebraic domains and guard sets, and parametric uncertainties. 

\exmpl{In this first example, we analyze Zeno stability of a hybrid system without time-invariant parametric uncertainties.
\\\noindent Consider the hybrid system $H = (Q,E,D,F,G,R)$, where 
\begin{itemize}
\item $Q = \{1,2,3\}$
\item $E = \{(1,2),(2,3),(3,1)\}$
\item  $D := \{D_1,D_2,D_3\}$ where
\aligneq{
&D_1 = \{x\in\Realsn{2}: x_1>0, x_2 + \frac{1}{2}x_1 \geq  0 \}\\
&D_2 = \{x\in\Realsn{2}: x_2-\frac{1}{2}x_1\geq 0, x_2 + \frac{1}{2}x_1 <  0 \}\\
&D_3 = \{x\in\Realsn{2}: x_1<0, x_2 + \frac{1}{2}x_1 \geq  0 \};}
\item $F = \{f_1,f_2,f_3\}$, where
\aligneq{
\dot{x} &= f_1(x) =\left (x_2,\;-5x_1-x\right )^T\\
\dot{x} &= f_2(x) = \left (-x_1^2-3,\;2x_2^2-\frac{1}{2}x_1^2\right )\\
\dot{x} &= f_3(x) = \left (x_2^2+x_1.\;-3x_1\right );
}
\item  $G := \{G_{12},G_{23}, G_{31}\}$ where 
\aligneq{
&G_{12}:=\left\{x\in \Realsn{2}: x_2\leq 0, \frac{1}{2}x_1+x_2=0 \right\}\\
&G_{23}:=\left\{x\in \Realsn{2}: x_2\leq 0, \frac{1}{2}x_1-x_2=0 \right\}\\
&G_{31}:=\left\{x\in \Realsn{2}: x_2> 0, x_1=0 \right\};
}
\item $R= \{\phi_{12}(x),\phi_{23}(x),\phi_{31}(x)\}$ where each $\phi_{ij}(x) = x$.
\end{itemize}
We note that this hybrid system is cyclic, as the pair $(Q,E)$ forms a directed cycle, with vertices $Q$ and edges $E$. A phase portrait of the system is given below in Figure~\ref{fig:nonlinear2}.
\begin{figure}[h]
\centering
\includegraphics[width=0.85\columnwidth]{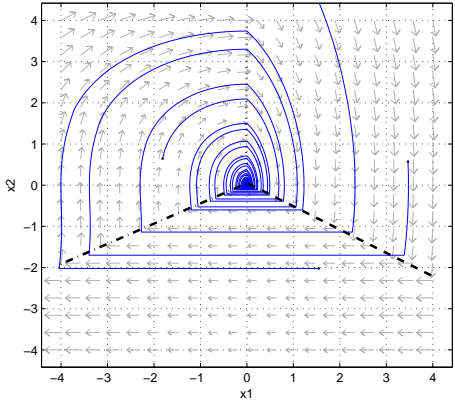}
\caption{Hybrid System in Example 1. Dashed line indicates $G_{12}$, dash-dotted line indicates $G_{23}$ and dotted line indicates $G_{31}$ }
\label{fig:nonlinear2}

\end{figure}

\noindent {\bf Results:}\\
 We wish to analyze Zeno stability for $z = \{z_1,z_2,z_3\}$, where $z_1 = z_2 = z_3 = [0,0]^T$. To solve Feasibility Problem 1,  we consider each 
\[W_q:= \mathbb{B}^2\cap D_q\]  
where $\mathbb{B}^2:=\{x\in \Realsn{2}: |x|\leq 1\}$.\\\\
We then search for 3 degree 8 polynomials to solve Feasibility Problem 1. Since we are able to solve Feasibility Problem 1 with such polynomials, we show using Theorem~\ref{robust} that $z = [0,0]^T$ is Zeno stable for $H$.
}

\exmpl{
Consider the hybrid system $H=(Q,E,D,F,G,R)$ with uncertain parameter $p\in (C,\infty)$ where
\begin{itemize}
\item $Q = \{1,2\}$
\item $E = \{(1,2),(2,1)$
\item $D = \{D_1,D_2\}$ where
\aligneq{
&D_1 := \{x\in \Realsn{2}: x_1+x_2\geq 0,px_1-x_2\geq 0\}\\
&D_2 := \Realsn{2}\backslash D_1
}
\item $F = \{f_1,f_2\}$ where
\aligneq{
&f_1 = \rmatrix{c}{-0.1\\
2}\\
&f_2 = \rmatrix{c}{- x_2 - x_1^3\\
x_1}}
\item $G = \{G_{12},G_{21}\}$ where
\aligneq{
&G_{12}:= \{x\in\Realsn{2}:x_2-px_1=0\}\\
&G_{21} := \{x\in \Realsn{2}: x_1+x_2 = 0\}
}
\item $R= \{\phi_{12}(x),\phi_{21}(x)\}$ where each $\phi_{ij}(x) = x$.
\end{itemize}

In this example, the uncertain parameter affects the switching rule. Provided below are simulations with 3 different fixed values of $p$. First, we consider the case when $p=1$, in Figure~\ref{fig:zenoswitch1}.
\begin{figure}[h]
\centering
\includegraphics[width=0.85\columnwidth]{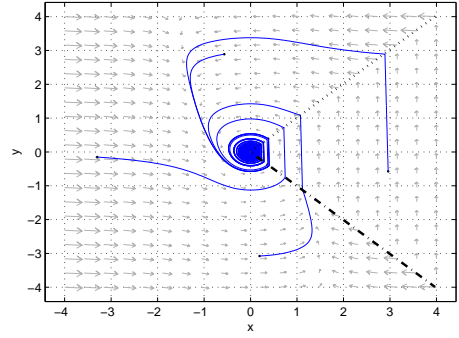}
\caption{Trajectories of Hybrid System in Example 2 with p=1. Dotted line indicates $G_{12}$ and dash-dotted line indicates $G_{21}$}

\label{fig:zenoswitch1}
\end{figure}

We see from inspection that the origin is Zeno stable. Furthermore, if we consider Figure~\ref{fig:zenoswitch3}, we see that even if we increase $p$ (thereby increasing the slope of $G_{21}$), we notice that the system remains Zeno stable.

\begin{figure}[h]
\centering
\includegraphics[width=0.85\columnwidth]{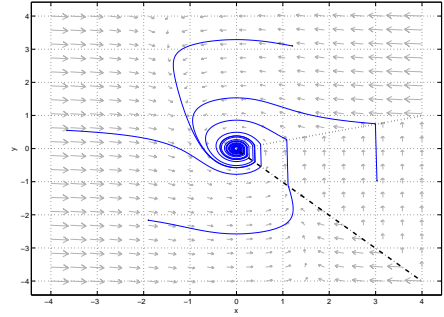}
\caption{Trajectories of Hybrid System in Example 2 with p=4. Dotted line indicates $G_{12}$ and dash-dotted line indicates $G_{21}$}

\label{fig:zenoswitch3}
\end{figure}

However, when we reduce the value of $p$, we notice that the system exhibits different asymptotic behavior. First, we see that if $p\leq -0.1$, the system will no longer exhibit Zeno stability. Indeed, in that circumstance, the system would not display any form of stable behavior. This is because the trajectories in $D_1$ would never reach the guard set (since the direction of the vector field would be parallel to the guard set). But even if $p\in(-0.1,1)$, we notice that the system asymptotically converges to limit cycles, as seen in Figure~\ref{fig:zenoswitch2}.

\begin{figure}[h]
\centering
\includegraphics[width=0.85\columnwidth]{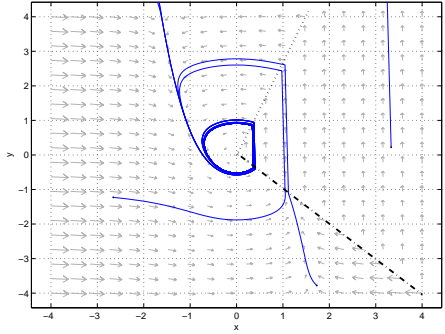}
\caption{Trajectories of Hybrid System in Example 2 with p=0.4. Dotted line indicates $G_{12}$ and dash-dotted line indicates $G_{21}$}

\label{fig:zenoswitch2}
\end{figure}

{\bf Results:}\\
We wish to analyze Zeno stability of $z = \{z_1,z_2\} = (0,0)$.
For our computational analysis, we first divide $D_2$ into $D_{21}$ and $D_{22}$,where
\aligneq{
&D_{21} := \{x\in \Realsn{2}: x_1+x_2\geq 0,-px_1+x_2\geq 0\}\\
&D_{22} := \{x\in \Realsn{2}: x_1+x_2\leq 0\}.
}
 We then search for a common Lyapunov function for both $D_{21}$ and $D_{22}$. The set of uncertain parameters is given by the inequality $P:=\{p\in\Reals: \tilde{p}:=p-C>0\}$, where $C$ is determined a priori. The goal is to find a lower bound on $C$ such that $z$ is Zeno stable. We use $W = W_1\cup W_2 = \{x\in\Realsn{2}: x_1^2+x_2^2< 5\}$. We then search for Lyapunov functions of varying degrees for different values of $C$. We note that as we increase the degree of $V_1$ and $V_{2}$, we are able to obtain a tighter lower bound on $C$. These results are given below in table~\ref{tab:tab1}.

\begin{table}[h]
\centering
\begin{tabular}{|c|c|}\hline
{\bf Degree of $V_1,V_2$} &{\bf Lower bound on $C$}\\\hline
8 & 2.11\\
10 & 1.87 \\
12 & 1.73\\\hline
\end{tabular}
\caption{Lower bound on $C$ for which $z$ is Zeno stable obtained for different degrees of $V_1$ and $V_2$}
\label{tab:tab1}
\vspace{-2mm}
\end{table}
We were unable to find a feasible $V_1$ and $V_2$ of degree less than 8. Unfortunately, we were unable to search for polynomials of degree greater than 12 owing to computational limitations.
}
\section{Conclusions}
In this paper, we present a Lyapunov based method for determining Zeno stability of compact sets in hybrid dynamical systems. The method presented makes use of the sum-of-squares decomposition, thus enabling the construction of higher-order Lyapunov functions. As such, the theorem presented can be used to certify Zeno stability in systems with nonlinear vector fields and reset maps, and time-invariant parametric uncertainties. The result can easily be simplified to certify Zeno stability for hybrid systems without parametric uncertainties as well. Examples of hybrid systems with polynomial domains, guard sets, vector fields and reset maps illustrating the use of the proposed method are also provided. We also provide an example of a nonlinear polynomial hybrid system with an uncertain parameter in the guard set. 

\bibliographystyle{IEEEtran}
\bibliography{mybib1}
\end{document}